\documentclass[12pt,reqno,a4]{amsart}
\usepackage{a4}
\usepackage{graphicx}
\usepackage{amsmath,amssymb,natbib,color,bbm,ifthen,graphicx,epsf}
\usepackage[latin1]{inputenc}

\def\Nset{{\mathbb{N}}}
\def\Rset{\mathbb R}
\def\Cset{\mathbb{C}}

\newenvironment{enum_A}
  {%
  \setlength{\leftmargini}{2em}\begin{enumerate}}
  {\end{enumerate}}

\newenvironment{enum_AP}
  {%
  \setlength{\leftmargini}{2em}\begin{enumerate}}
  {\end{enumerate}}

\def\TV{\mathrm{TV}}
\def\1{\mathbbm{1}}
\newcommand{\PE}{\mathbb{E}}
\def\PP{\mathbb P}

\def\re{\mathrm{Re}}

\def\rme{\mathrm{e}}
\def\ie{\textit{i.e.}}
\newcommand{\wrt}{with respect to}
\newcommand{\as}{a.s.}
\newcommand{\Pdouble}{\check{P}}
\newcommand{\PPdouble}{\check{\PP}}
\newcommand{\PEdouble}{\check{\PE}}
\newcommand\driftfunction{W}

\def\mcf{\mathcal{F}}

\def\Xsigma{\mathcal{X}}

\def\eqsp{\;}
\newcommand{\eqdef}{\ensuremath{\stackrel{\mathrm{def}}{=}}}
\newcommand{\driftU}{U}
\newcommand{\driftV}{V}
\newcommand{\driftv}{v}
\newcommand{\rphiunnorm}{\phi \circ H_\phi^{-1}}
\newcommand{\rphinorm}{r_\phi}

\def\Xset{\mathsf{X}}

\newtheorem{theo}{Theorem}[section]

\newtheorem{prop}[theo]{Proposition}

 \theoremstyle{remark}
\newtheorem{rem}{Remark}

\numberwithin{equation}{section}

\date{\today\ (draft version)}

\title[Computable Convergence Rates for Markov Chains]{Computable Convergence Rates for Subgeometrically Ergodic Markov Chains}

\author{Randal Douc$^\star$}
\address{Randal DOUC (Corresponding Author), CMAP, Ecole Polytechnique, 91128 Palaiseau Cedex, France}
\email{douc@cmapx.polytechnique.fr}
\author{Eric Moulines}
\address{Eric MOULINES, Ecole Nationale Supérieure des Télécommunications, Département de Traitement du Signal et des Images, 46, rue Barrault, 75634 PARIS Cédex 13, France}
\email{moulines@tsi.enst.fr}
\author{Philippe Soulier}
\address{Philippe SOULIER, Department of Mathematics of University Paris X, Bâtiment G, Bureau E15, 200 avenue de la République
92000 Nanterre Cedex, France}
\email{philippe.soulier@u-paris10.fr}

\begin{document}

\maketitle

\begin{abstract}
  In this paper, we give quantitative bounds on the $f$-total variation distance
  from convergence of an Harris recurrent Markov chain on an arbitrary under
  drift and minorisation conditions implying ergodicity at a sub-geometric rate.
  These bounds are then specialized to the stochastically monotone case, covering the case where
  there is no minimal reachable element. The results are
  illustrated on two examples from queueing theory and Markov Chain Monte Carlo.
\end{abstract}



 \section{Introduction}
 Let $P$ be a Markov transition kernel on a state space $\Xset$
 equipped with a countably generated $\sigma$-field $\Xsigma$.
 For a control function $f: \Xset \to [1,\infty)$, the $f$-\textit{total variation} or $f$-\textit{norm} of a
 signed measure $\mu$ on $\Xsigma$ is defined as
\[
\| \mu \|_f := \sup_{|g| \leq f} | \mu(g) | \eqsp.
\]
When $f \equiv 1$, the $f$-norm is the total variation norm, which is denoted
$\| \mu \|_{\TV}$.
%
%
We assume that $P$ is aperiodic positive Harris recurrent with stationary distribution $\pi$. Our goal is to obtain \emph{quantitative}
bounds on convergence rates, \ie\ rate of the form
\begin{equation}  \label{eq:protoresult}
r(n) \| P^n(x,\cdot) - \pi \|_f \leq g(x) \eqsp, \quad \text{for all $x \in \Xset$}
\end{equation}
where  $f$ is a control function $f: \Xset \to [1,\infty)$, $\{ r(n) \}_{n \geq 0}$ is a non-decreasing sequence,
and $g: \Xset  \to [0,\infty]$ is a function which can be computed explicitly.
As emphasized in \cite[section 3.5]{roberts:rosenthal:2004}, quantitative bounds have a substantial
history  in Markov chain theory. Applications are numerous including convergence analysis of Markov Chain Monte Carlo (MCMC) methods,
transient analysis of queueing systems or storage models, etc.
With few exception however, these quantitative bounds were derived under conditions which imply
geometric convergence, \textit{i.e.} $ r(n) = \beta^n$, $n \geq 1$.
(see for instance \cite{meyn:tweedie:1994}, \cite{rosenthal:1995},
\cite{roberts:tweedie:1999}, \cite{roberts:rosenthal:2004}, and \cite{baxendale:2005}).

In this paper, we study conditions under which \eqref{eq:protoresult} hold for
sequences in the set $\Lambda$ of subgeometric rate functions from
\cite{nummelin:tuominen:1983}, defined as the family of sequences $\{ r(n) \}_{n \geq 0}$ such that
$r(n)$ is non decreasing and $\log r(n)/n \downarrow 0$ as $n \to \infty$.
Without loss of generality, we assume that $r(0)=1$ whenever $r \in \Lambda$.
These rates of convergence 
have been only scarcely considered in the literature.
Let us briefly summarize the results available for
convergence at subgeometric rate for general state-space chain.  To our best knowledge, the first
result for subgeometric sequence has been obtained by \cite{nummelin:tuominen:1983},
who derive sufficient conditions for  $\| \xi P^n - \pi \|_{\TV}$ to be of order $o(r^{-1}(n))$.
The basic condition involved in this work is the ergodicity of order $r$ (or $r$-ergodicity), defined as
\begin{equation}\label{eq:r-moment}
 \sup_{x \in B} \PE_x \left[ \sum_{k=0}^{\tau_B-1} r(k) \right] < \infty \eqsp.
\end{equation}
where $\tau_B \eqdef \inf \{ n \geq 1, X_k \in B \}$ (with the convention that $\inf \emptyset = \infty$)
is the return time to some accessible some small set $B$ (\ie\ $\pi(B) > 0$).
These results were later extended by  \cite{tuominen:tweedie:1994} to $f$-norm for general control
functions $f : \Xset \to [1,\infty)$ under $(f,r)$-ergodicity, which states that
\begin{gather} \label{eq:fr-moment}
  \sup_{x \in B} \PE_x \left [ \sum_{k=0}^{\tau_B-1} r(k) f(X_k) \right] < \infty
\end{gather}
for some accessible small set $B$. These contributions do not provide computable expressions for the bounds in \eqref{eq:protoresult}.

A direct route to quantitative bounds for subgeometric sequences has been opened by
\cite{veretennikov:1997,veretennikov:1999}, based on  coupling techniques
(see \cite{gulinsky:veretennikov:1993} and \cite{rosenthal:1995} for the coupling construction of Harris
recurrent Markov chains).
This method consists in relating the bounds \eqref{eq:protoresult} to
a moment of the \textit{coupling time} through Lindvall's
inequality \cite{lindvall:1979,lindvall:1992}.
\cite{veretennikov:1997,veretennikov:1999} focus on a particular class
of Markov chains, the so-called functional autoregressive processes,
defined as $X_{n+1} = g(X_n) + W_{n+1}$, where $g : \Rset^d \to
\Rset^d$ is a Borel function and $(W_n)_{n \geq 0}$ is an i.i.d.
sequence, and provides expressions of the bounds in
\eqref{eq:protoresult} with the total variation distance ($f \equiv
1$) and polynomial rate functions $r(n)= n^\beta$, $n \geq 1$.  These
results have later been extended, using similar techniques, to
\textit{truly} subgeometric sequence, \textit{i.e.}  $\{ r(n) \}_{n \geq 0} \in \Lambda$ satisfying $\lim_{n \to \infty} r(n)
n^{-\kappa} = \infty$ for any $\kappa$, in
\cite{klokov:veretennikov:2004}, for a more general class of
functional autoregressive process.

\cite{fort:moulines:2003:SPA} derived quantitative bounds of the form
\eqref{eq:protoresult} for possibly unbounded control functions and
polynomial rate functions, also using the coupling method.
The bound for the modulated moment of the coupling time is obtained
from a particular drift condition introduced by \cite{fort:moulines:2000}
later extended by \cite{jarner:roberts:2001}.  This method is based
on a recursive computation of the polynomial moment of the coupling time (see \cite[proposition 7]{fort:moulines:2003})
which is related to the moments of the hitting time of a bivariate chain to a set where coupling might occur.
This proof is tailored to the polynomial case and cannot be easily adapted
to the general subgeometric case (see \cite{fort:2001} for comments).

The objective of this paper is to generalize the results mentioned
above in two directions.  We consider Markov chains over general
state space  and we study general subgeometrical rates of
convergence instead of polynomial rates
\cite{fort:moulines:2003:SPA}. We establish  a  family of
convergence bound (with a trade-off between the rate and the norm)
extending to the subgeometrical case the computable bounds obtained
in the geometrical case by \cite{rosenthal:1995} and later refined
by \cite{roberts:tweedie:1999} and \cite{douc:moulines:rosenthal:2004} (see \cite[Theorem 12]{roberts:rosenthal:2004} and the references therein).
The method, based on  coupling associated,  provides a short and nearly self-contained proof of
the results presented in \cite{nummelin:tuominen:1983} and \cite{tuominen:tweedie:1994}: this allows for
intuitive understanding of these results, while also avoiding various analytic technicalities of the previous proofs of these
theorems.

The paper is organized as follows. In section \ref{sec:statementoftheresults},
we present our assumptions and  state our main results. In section \ref{sec:stochastically:monotone}, we specialize our result to
stochastically monotone Markov chains and derive bounds which extends results reported earlier by
\cite{scott:tweedie:1996} and \cite{roberts:tweedie:2000}. Examples from queueing theory and MCMC
are discussed in section \ref{sec:applications} to support our findings and illustrate the numerical computations of
the bounds.

\section{Statements of the results}
\label{sec:statementoftheresults}
The proof is based on the coupling construction (briefly recalled in section \ref{sec:proofmainresult}). It is assumed that the chain admits a small set:
\begin{enum_A}
\item \label{assumption:smallset} There exist a set $C \in \Xsigma$, a constant $\epsilon > 0$ and a probability measure $\nu$ such that,
for all $x \in C$, $P(x,\cdot) \geq \epsilon \nu(\cdot)$.
\end{enum_A}
For simplicity, only one-step minorisation is considered in this paper.
Adaptations to $m$-step minorisation can be carried out as in
\cite{rosenthal:1995} (see also \cite{fort:2001} and \cite{fort:moulines:2003:SPA}).

Let $\Pdouble$ be a Markov transition kernel on $\Xset \times \Xset$ such that, for all $A \in \Xsigma$,
\begin{align}
\label{eq:defPdouble}
&\Pdouble(x,x', A \times \Xset)= P(x,A) \1_{(C \times C)^c}(x,x') + Q(x,A) \1_{C \times C}(x,x') \\
&\Pdouble(x,x', \Xset \times A)= P(x', A) \1_{(C \times C)^c}(x,x') + Q(x',A) \1_{C \times C}(x,x')
\end{align}
where $A^c$ denotes the complementary of the subset $A$ and $Q$ is the so-called residual kernel defined, for $x \in C$ and $A \in \Xsigma$ by
\begin{equation}
\label{eq:DefinitionResidualKernel}
Q(x,A) = \begin{cases}
(1-\epsilon)^{-1} \left( P(x,A) - \epsilon \nu(A) \right) & 0 < \epsilon < 1 \\
\nu(A) & \epsilon = 1
\end{cases}
\end{equation}
One may for example set
\begin{multline}
\label{eq:bivariatekernel-independentcomponent}
\Pdouble(x,x'; A \times A') = \\
P(x,A) P(x',A') \1_{(C \times C)^c}(x,x') + Q(x,A) Q (x',A) \1_{C \times C}(x,x') \eqsp,
\end{multline}
but, as seen below, this choice is not always the most suitable.
For $(x,x') \in \Xset \times \Xset$,
denote by $\PPdouble_{x,x'}$ and $\PEdouble_{x,x'}$ the law and the expectation of a Markov chain with initial
distribution $\delta_x \otimes \delta_{x'}$ and transition kernel $\Pdouble$.

Our second condition is a bound on the  moment of the hitting time of the bivariate chain to  $C \times C$ under the probability $\PPdouble_{x,x'}$.
Let $\{ r(n) \} \in \Lambda$ be a subgeometric sequence and set: $R(n) \eqdef \sum_{k=0}^{n-1} r(k)$.
Denote by $\sigma_{C \times C} \eqdef \inf \{ n \geq 0, (X_n,X'_n) \in {C \times C} \}$ the first hitting time of  ${C \times C}$ and let
\begin{equation}
\label{eq:definition:driftr}
\driftU(x,x') \eqdef \PEdouble_{x,x'} \left[ \sum_{k=0}^{\sigma_{C \times C}} r(k)  \right] \eqsp.
\end{equation}
Let  $\driftv: \Xset \times \Xset \to [0,\infty)$ be a measurable function and set
\begin{equation}
\label{eq:definition:driftw}
\driftV(x,x') = \PEdouble_{x,x'} \left[ \sum_{k=0}^{\sigma_{C \times C}} \driftv(X_k,X'_k)\right]
\end{equation}
\begin{enum_A}
\setcounter{enumi}{1}
\item
\label{assumption:rate}
For any $(x,x') \in \Xset \times \Xset$, $\driftU(x,x') < \infty$ and
\begin{equation}
\label{eq:assumption:rate}
b_{\driftU} \eqdef \sup_{(x,x') \in C \times C} \Pdouble \driftU(x,x') = \sup_{(x,x') \in {C \times C}} \PEdouble_{x,x'} \left[ \sum_{k=0}^{\tau_{C \times C}-1} r(k) \right] < \infty
\end{equation}
\item \label{assumption:modulatedmoment}
For any $(x,x') \in \Xset \times \Xset$, $\driftV(x,x') < \infty$ and
\begin{equation}
\label{eq:assumption:modulatedmoment}
b_{\driftV} = \sup_{(x,x') \in {C \times C}} \Pdouble \driftV(x,x')= \sup_{(x,x') \in {C \times C}} \PEdouble_{x,x'} \left[ \sum_{k=1}^{\tau_{C \times C}} \driftv(X_k,X'_k) \right]  < \infty \eqsp.
\end{equation}
\end{enum_A}

%
We will establish that $R$ is the maximal rate of convergence
(that can be deduced from assumptions \ref{assumption:smallset}-\ref{assumption:modulatedmoment}) and that this rate
is associated to convergence in total variation norm. On the other hand, we will show that the difference $P(x,\cdot) - P(x',\cdot)$ remains bounded
in $f$-norm for any function $f$  satisfying $f(x) + f(x') \leq \driftV(x,x')$ for any $(x,x') \in \Xset \times \Xset$.
Using an interpolation technique, we will derive rate of convergence
$1 \leq s \leq r$ associated to some $g$-norm, $0 \leq g \leq f$. To construct such interpolation, we consider pair of positive functions
$(\alpha,\beta)$ satisfying, for some $0 \leq \rho \leq 1$,
\begin{equation}
\label{eq:additivity}
\alpha(u) \beta(v) \leq \rho u+ (1-\rho) v \eqsp, \quad \text{for all $(u,v) \in \Rset^+ \times \Rset^+$} \eqsp.
\end{equation}
Functions satisfying this condition can be obtained from Young's inequality. Let $\psi$ be a real valued, continuous, strictly increasing function on $\Rset^+$
such that $\psi(0)= 0$; then for any $(a,b) > 0$,
$$
ab \leq \mathcal{P}(a) + \mathcal{F}(b) \eqsp, \text{where} \quad \mathcal{P}(a) \eqdef \int_0^a \psi(x) dx \quad \text{and} \quad \mathcal{F}(b) \eqdef \int_0^b \psi^{-1}(x) dx \eqsp,
$$
where $\psi^{-1}$ is the inverse function of $\psi$. If we set $\alpha(u) \eqdef \mathcal{P}^{-1}(\rho u)$ and $\beta(v)= \mathcal{F}^{-1}((1-\rho)v)$, then the pair $(\beta,\alpha)$
satisfies \eqref{eq:additivity}. Taking $\psi(x)= x^{p-1}$ for some $p \geq 1$ gives the special case $\left\{ (p \rho u)^{1/p}, (p (1- \rho) u/(p-1))^{(p-1)/p} \right\}$.
\begin{theo} \label{theo:mainresult}
Assume \ref{assumption:smallset}, \ref{assumption:rate}, and \ref{assumption:modulatedmoment}.
Define
\begin{equation}
\label{eq:constantM}
M_\driftU \eqdef \sup_{k \in \Nset} \left\{ \left( b_\driftU r(k) \frac{1-\epsilon}{\epsilon} -R(k+1) \right)_+ \right\}  \quad \text{and} \quad
M_\driftV \eqdef  b_{\driftV} \frac{1-\epsilon}{\epsilon} \eqsp,
\end{equation}
where $(x)_+ \eqdef \max(x,0)$. Then, for any $(x,x') \in \Xset \times \Xset$,
\begin{align}
& \| P^n(x,\cdot) - P^n(x',\cdot) \|_{\TV} \leq \frac{\driftU(x,x')+ M_{\driftU}}{R(n) + M_\driftU} \\
& \| P^n(x,\cdot) - P^n(x',\cdot) \|_{f} \leq \driftV(x,x') + M_{\driftV} \eqsp,
\end{align}
for any non-negative function $f$ satisfying, for any $(x,x') \in \Xset \times \Xset$,  $f(x) + f(x') \leq \driftV(x,x') + M_\driftV$.
Let $(\alpha,\beta)$ be two positive functions satisfying \eqref{eq:additivity} for some $0 \leq  \rho \leq 1$.
Then, for any  $(x,x') \in \Xset \times \Xset$ and  $n \geq 1$, :
  \begin{align} \label{eq:interpolationfnorm}
    \| P^n(x,\cdot) - P^n(x',\cdot) \|_{g} \leq \frac{\rho \left(\driftU(x,x') + M_{\driftU}\right) + (1-\rho) \left( \driftV(x,x') + M_{\driftV} \right) }{\alpha \circ \left\{ R(n) + M_{\driftU} \right\}}
  \end{align}
for any non-negative function $g$ satisfying, for any $(x,x') \in \Xset \times \Xset$, $g(x) + g(x') \leq \beta \circ \left\{ \driftV(x,x') + M_\driftV \right\}$.
\end{theo}
The proof is postponed to section \ref{sec:proofmainresult}.
\begin{rem}
Because the sequence $\{ r(k) \}$ is subgeometric, $\lim_{k \to \infty} r(k)/R(k+1)=0$. Therefore, the sequence $\{b_U r(k) (1-\epsilon)/\epsilon - R(k)\}$
has only finitely many non-negative terms, which implies that $M_U < \infty$.
\end{rem}
\begin{rem}
When assumption \ref{assumption:rate}, then \ref{assumption:modulatedmoment} is automatically  satisfied for some function $\driftv$. Note that
\[
\PEdouble_{x,x'} \left[ \sum_{k=0}^{\sigma_{C \times C}} r(k) \right] = \PEdouble_{x,x'} \left[ \sum_{k=0}^{\sigma_{C \times C}} r(\sigma_{C \times C}-k) \right] \eqsp.
\]
On the other hand, for all $(x,x') \in \Xset \times \Xset$,
\begin{multline*}
\PEdouble_{x,x'} \left[ r(\sigma_{C \times C}-k) \1_{\{\sigma_{C \times C} \geq k\}} \right]
\\ = \PEdouble_{x,x'} \left[ \PEdouble_{X_k,X'_k} \left[ r(\sigma_{C \times C}) \right] \1_{ \{ \sigma_{C \times C} \geq k \}} \right] =  \PEdouble_{x,x'} \left[ \driftv_{r}(X_k,X'_k) \1_{\{\sigma_{C \times C} \geq  k\}} \right] \eqsp,
\end{multline*}
where $\driftv_r(x,x') \eqdef \PEdouble_{x,x'} [ r(\sigma_{C \times C})] $. This relation implies that
$$\PEdouble_{x,x'} \left[ \sum_{k=0}^{\sigma_{C \times C}} r(k) \right]= \PEdouble_{x,x'} \left[ \sum_{k=0}^{\sigma_{C \times C}} \driftv_r(X_k,X'_k) \right] \eqsp,
\quad \text{for all $(x,x') \in \Xset \times \Xset$} \eqsp.
$$
However, in particular when using drift functions, it is sometimes easier to apply theorem \ref{theo:mainresult} with function a function $\driftv$ which does not coincide
with $\driftv_r$.
\end{rem}

%

To check assumptions \ref{assumption:rate} and \ref{assumption:modulatedmoment} it is often useful to use a drift conditions.
Drift conditions implying convergence at polynomial rates have been recently proposed in \cite{jarner:roberts:2001}. These conditions
have later been extended to general subgeometrical rates by \cite{douc:fort:moulines:soulier:2004}.
Define by $\mathcal{C}$ the set of functions
\begin{multline}
\label{eq:definitionsetC}
\mathcal{C} \eqdef \left\{ \phi: [1,\infty) \to \Rset^+ \eqsp, \text{$\phi$ is concave, differentiable and} \right.\\
\left.   \phi(1) > 0, \lim_{v\to\infty} \phi(v) = \infty, \lim_{v\to\infty}\phi'(v) = 0  \right\} \eqsp.
\end{multline}
For $\phi \in \mathcal{C}$, define  the function $H_\phi: [1,\infty) \to [0,\infty)$ as $H_\phi (v) \eqdef \int_1^v \frac{dx}{\phi(x)}$.
 Since $\phi$ is non decreasing, $H_\phi$ is a non decreasing concave differentiable function on $[1, \infty)$ and
 $\lim_{v \to \infty} H_\phi(v)= \infty$. The inverse  $H_\phi^{-1}: [0, \infty) \to [1, \infty)$ is also an
 increasing and differentiable function, with derivative $(H_\phi^{-1})' = \phi\circ H_\phi^{-1}$.
Note that
$(\log\{ \phi\circ H_\phi^{-1} \})' = \phi'\circ H_\phi^{-1}$. Since $H_\phi$ is increasing and $\phi'$ is
 decreasing,  $\rphiunnorm$ is log-concave, which implies that the sequence
\begin{equation}
\label{eq:definitionphi}
 \rphinorm(n) \eqdef  \rphiunnorm(n) / \rphiunnorm(0) \eqsp,
\end{equation}
belongs to the set of subgeometric sequences $\Lambda$. Consider the following assumption
\begin{enum_A}
\setcounter{enumi}{3}
\item \label{assumption:driftdouble}
There exists a function $\driftfunction: \Xset \times \Xset \to [1,\infty)$, a function $\phi \in \mathcal{C}$ and a constant
$b$ such that $\Pdouble \driftfunction(x,x') \leq \driftfunction(x,x') - \phi \circ \driftfunction(x,x')$ for $(x,x') \not \in C \times C$
and $\sup_{(x,x') \in C \times C} \Pdouble \driftfunction(x,x') < \infty$.
\end{enum_A}
It is shown in \cite{douc:fort:moulines:soulier:2004} that under \ref{assumption:driftdouble}, \ref{assumption:rate} and
\ref{assumption:modulatedmoment} are satisfied with the rate sequence $r_\phi$ and the control function $\driftv = \phi \circ \driftfunction$.
In addition, it is possible to deduce explicit bounds for the constants $B_{\driftU}$, $b_{\driftU}$, $B_{\driftV}$ and $b_{\driftV}$ from the
constants appearing in the drift condition.
\begin{prop}
\label{prop:mainresultwithdrift}
Assume \ref{assumption:driftdouble}. Then, \ref{assumption:rate} and \ref{assumption:modulatedmoment} hold with $\driftv = \phi \circ \driftfunction$,
$r = r_\phi$ and
\begin{align}
\label{eq:boundV0}
&\driftU(x,x') \leq 1  + \frac{\rphinorm(1)}{\phi(1)} \left\{ \driftfunction(x,x') - 1 \right\} \1_{(C \times C)^c}(x,x') \eqsp, \\
\label{eq:boundV1}
&\driftV(x,x') \leq \sup_{{C \times C}} \phi \circ \driftfunction + \driftfunction(x,x') \1_{(C \times C)^c}(x,x') \eqsp, \\
\label{eq:boundV0-1}
&b_{\driftU} \leq 1  + \frac{\rphinorm(1)}{\phi(1)} \left\{ \sup_{C \times C} \Pdouble \driftfunction - 1 \right\} \eqsp,\\
&b_{\driftV} \leq \sup_{C \times C} \phi \circ \driftfunction   +  \sup_{C \times C} \Pdouble \driftfunction \eqsp.
\end{align}
\end{prop}
The proof is in section \ref{sec:proofmainresults-auxiliary}. Proposition \ref{prop:mainresultwithdrift} is only
partially satisfactory because Assumption \ref{assumption:driftdouble} is formulated on the bivariate kernel $\Pdouble$.
It is in general easier to establish directly the drift condition on the kernel $P$ and to deduce from this condition a
drift condition for an appropriately defined kernel $\Pdouble$ (see \cite[Proposition 11]{roberts:rosenthal:2004} for a similar construction for
geometrically ergodic Markov chain). Consider the following assumption:
\begin{enum_A}
\setcounter{enumi}{4}
\item \label{assumption:driftsimple}
There exists a  function $\driftfunction_0: \Xset \times \Xset \to [1,\infty)$, a function $\phi_0 \in \mathcal{C}$ and a constant
$b_0$ such that $P \driftfunction_0 \leq \driftfunction_0 - \phi_0 \circ \driftfunction_0 + b_0 \1_{C}$.
\end{enum_A}

\begin{theo}
\label{theo:singledriftcondition}
Suppose that \ref{assumption:smallset} and \ref{assumption:driftsimple} are satisfied. Let $d_0 \eqdef \inf_{x \not \in C} \driftfunction_0(x)$.
Then, if $\phi_0(d_0) > b_0$, the kernel $\Pdouble$ defined in \eqref{eq:bivariatekernel-independentcomponent} satisfies the bivariate drift condition
\ref{assumption:driftdouble} with
\begin{align}
&\driftfunction(x,x')= \driftfunction_0(x) + \driftfunction_0(x') -1 \\
&\phi = \lambda \phi_0 \eqsp, \quad \text{for any $\lambda \eqsp, 0 < \lambda < 1 - b_0 / \phi_0(d_0) $} \\
&\sup_{C \times C} \Pdouble \driftfunction \leq 2 (1-\epsilon)^{-1} \left\{  \sup_C P \driftfunction_0 - \epsilon \nu(\driftfunction_0) \right\}  -1 \eqsp.
\end{align}
where the kernel $Q$ is defined in \eqref{eq:DefinitionResidualKernel}.
\end{theo}
The proof is postponed to the appendix.
\begin{rem}
  Since the function $\phi_0$ is non-decreasing and $\lim_{v \to \infty}
  \phi_0(v) = \infty$, one may always find $d$ such that the condition
  $\phi_0(1)+\phi_0(d)\geq b_0(1-\alpha)^{-1}+2$ is fulfilled. The assumptions
  of the theorem above are satisfied provided that the associated level set $\{
  V_0 \leq d \}$ is small. This will happen of course if all the level sets are
  $1$-small, which may appear to be a rather strong requirement. More realistic
  conditions may be obtained by using small sets associated to the iterate $P^m$ of the kernel
  (see \textit{e.g.} \cite{rosenthal:1995},
  \cite{fort:2001} and \cite{fort:moulines:2003:SPA}).
\end{rem}

\subsection{Stochastically ordered chains}
\label{sec:stochastically:monotone}
In this section, we show how to define the kernel $\Pdouble$ and obtain a drift
condition for stochastically ordered Markov chain.  Let $\Xset$ be a totally ordered set, and denote $\preceq$ the order
relation. For $a\in\Xset$, denote $(-\infty,a]=\{x\in\Xset:\, x\preceq
a\}$ and $[a,+\infty) = \{x\in\Xset:\, a \preceq x\}$.
A transition kernel $P$ on $\Xset$ is called \textit{stochastically
  monotone} if  for all $a \in \Xset$, $P(\cdot,(-\infty,a])$ is non increasing.
Stochastic monotonicity has been seen to be crucial in the analysis of queuing network, Markov Monte-Carlo methods, storage models, etc.
Stochastically ordered Markov chains have been considered in \cite{lund:tweedie:1996},
\cite{lund:meyn:tweedie:1996}, \cite{scott:tweedie:1996} and
\cite{roberts:tweedie:2000}.  In the first two papers, it is assumed that there
exists an atom at the bottom of the state space.
\cite{lund:meyn:tweedie:1996} cover only
geometric convergence; subgeometric rate of convergence are considered in
\cite{scott:tweedie:1996}.  \cite{roberts:tweedie:2000} covers the case where
the bottom of the space is a small set  but restrict their
attentions to conditions implying geometric rate of convergence.

For a general stochastically monotone Markov kernel $P$, it is always possible
to define the bivariate kernel $\Pdouble$ (see \eqref{eq:defPdouble})
so that the two components $\{X_n\}_{n \geq 0}$ and $\{X'_n\}_{n \geq 0}$
are \emph{pathwise ordered}, \ie\ their initial order is preserved at all times.

The construction goes as follows.
For $x \in \Xset$, $u \in [0,1]$ and $K$ a transition kernel on $\Xset$
denote by $G_K^-(x,u)$ the quantile function associated to the probability measure~$K(x,\cdot)$
\begin{gather} \label{eq:defquantile}
G_K^-(x,u) = \inf \{ y \in \Xset, K(x,(-\infty,y]) \geq u \} \eqsp.
\end{gather}
Assume that \ref{assumption:smallset} holds.
For $(x,x') \in \Xset \times \Xset$ and  $A \in \Xsigma \otimes \Xsigma$, define the transition kernel $\Pdouble$ by
\begin{multline*}
\1_{(C \times C)^c(x,x')} \Pdouble(x,x';A) = \int_0^1 \1_{A}(G_P^-(x,u),G_P^-(x',u)) \, du  \\
+ \1_{C \times C} (x,x') \; \int_0^1 \1_{A}(G_Q^-(x,u),G_Q^-(x',u)) \, du  \eqsp,
\end{multline*}
where $Q$ is the residual kernel defined in \eqref{eq:DefinitionResidualKernel}.
It is easily seen that, by construction, the set $\{(x,x')\in\Xset\times \Xset:\, x \preceq x'\}$ is absorbing for the
kernel $\Pdouble$.
%

In the sequel, we assume that \ref{assumption:smallset} holds for some $C \eqdef (-\infty,x_0]$
(\ie\ that there is a small set at the bottom of the space).
Let $\driftv_0: \Xset \to [1,\infty)$ be a measurable function and define:
\begin{equation}
\label{eq:SingleDriftSm}
\driftU_0(x) \eqdef \PE_x \left[ \sum_{k=0}^{\sigma_C} r(k) \right] \quad \text{and} \quad \driftV_0(x)= \PE_x \left[ \sum_{k=0}^{\sigma_C} \driftv_0(X_k) \right] \eqsp.
\end{equation}
Consider the following assumptions:
\begin{enum_AP}
\setcounter{enumi}{1}
\item \label{assumption:rate-stochasticallymonotone } For any $x \in \Xset$, $\driftU_0(x) < \infty$ and $\sup_{C} Q \driftU_0 = b_{\driftU_0} < \infty$,
\item \label{assumption:modulatedmoment-stochasticallymonotone}
For any $x \in \Xset$, $\driftV_0(x) < \infty$ and $\sup_{C} Q \driftV_0 = b_{\driftV_0} < \infty$,
\end{enum_AP}
\begin{theo}
\label{theo:mainresult:stochasticallymonotone}
Assume that \ref{assumption:smallset}-\ref{assumption:rate-stochasticallymonotone }-\ref{assumption:modulatedmoment-stochasticallymonotone}  holds for some set $C \eqdef (-\infty, x_0]$.
Then, \ref{assumption:rate} and \ref{assumption:modulatedmoment} hold with $\driftU(x,x') = \driftU_0(x \vee x')$, $\driftV(x,x') = \driftV_0(x \vee x')$,
$\driftv(x,x') = \driftv_0(x \vee x')$, $b_\driftU = b_{\driftU_0}$ and $b_\driftV = b_{\driftV_0}$.
\end{theo}
The proof is obvious and omitted for brevity. As mentioned above, drift conditions often provide an easy path to prove conditions such as \ref{assumption:rate-stochasticallymonotone } and
\ref{assumption:modulatedmoment-stochasticallymonotone}. Consider the following assumption:
\begin{enum_AP}
\setcounter{enumi}{3}
\item \label{assumption:driftfunction-stochasticallymonotone}
There exists a a nonnegative function $\driftfunction_0: \Xset \to [1,\infty)$, a function $\phi \in \mathcal{C}$
such that for $x \not \in C$, $P \driftfunction_0 \leq \driftfunction_0 - \phi \circ \driftfunction_0$ and $\sup_C P \driftfunction_0 < \infty$.
\end{enum_AP}
Using, as above \cite{douc:fort:moulines:soulier:2004}, it may be shown that this assumption implies \ref{assumption:rate-stochasticallymonotone } and
\ref{assumption:modulatedmoment-stochasticallymonotone} and allows to compute explicitly the constants.
\begin{theo}
\label{theo:mainresultwithdrift-stochasticallymonotone}
Assume \ref{assumption:smallset} and \ref{assumption:driftfunction-stochasticallymonotone}. Then \ref{assumption:rate-stochasticallymonotone } and \ref{assumption:modulatedmoment-stochasticallymonotone}
hold with $\driftv_0 = \phi \circ \driftfunction_0$, $r = r_\phi$, and
\begin{align}
\label{eq:boundV0-stochasticallymonotone}
&\driftU_0(x) \leq 1 + \frac{\rphinorm(1)}{\phi(1)} \left\{  \driftfunction_0(x) - 1 \right\} \1_{C^c}(x) \eqsp, \\
\label{eq:boundV1-stochasticallymonotone}
&\driftV_0(x) \leq \sup_{C} \phi \circ \driftfunction_0 + \driftfunction_0(x) \1_{C^c}(x) \eqsp, \\
\label{eq:boundV0-1-stochasticallymonotone}
&b_{\driftU_0} \leq 1 + \frac{\rphinorm(1)}{\phi(1)} \left((1-\epsilon)^{-1} \left\{\sup_{C} P \driftfunction_0-\epsilon \nu (\driftfunction_0) \right\}- 1  \right) \\
\label{eq:boundV1-1-stochasticallymonotone}
&b_{\driftV_0} \leq   \sup_{C}  \phi \circ \driftfunction_0+ (1-\epsilon)^{-1} \left\{\sup_{C} P \driftfunction_0-\epsilon \nu (\driftfunction_0) \right\}.
\end{align}
\end{theo}
The proof is entirely similar to Proposition \ref{prop:mainresultwithdrift} and is omitted.

\section{Applications}
\label{sec:applications}
\subsection{the embedded M/G/1 queue}
\label{subsec:MG1queue}
In a M/G/1  queue, customers arrive into a service operation according to a Poisson process with parameter $\lambda$.
Customers bring  jobs requiring a service times which are independent of each others and of the inter-arrival time with
common distribution $B$ concentrated on $(0,\infty)$ (we assume that the service time distribution has no probability mass at $0$).
Consider the random variable $X_n$ which counts customers immediately after each service time ends.
$\{ X_n \}_{n \geq 0}$ is a Markov chain on integers with transition matrix
\begin{equation}
\label{eq:M/G/1-transition-matrix}
P=
\left(
  \begin{array}{ccccc}
    a_0 & a_1 & a_2 & a_3 & \dots \\
    a_0 & a_1 & a_2 & a_3 & \dots \\
    0 & a_0 & a_1 & a_2 & \dots \\
    0 & 0 & a_0 & a_1 & \dots \\
    \vdots & \ddots& \ddots & \ddots & \ddots \\
  \end{array}
\right)
\end{equation}
where for each $j \geq 0$, $a_j \eqdef \int_0^\infty \left\{ \rme^{-\lambda t} (\lambda t)^j/ j!\right\} dB(t)$ (see \cite[proposition 3.3.2]{meyn:tweedie:1993}).
It is known that $P$ is irreducible, aperiodic, and positive recurrent if $\rho \eqdef \lambda m_1 = \sum_{j=1}^\infty j a_j < 1$, where for $u > 0$,
$m_u \eqdef \int t^u d B(t)$.
Applying the results derived above, we will compute explicit bounds  (depending on $\lambda$, $x$ and the moments of the service time distribution)
for the convergence bound $\| P^n(x,\cdot) - \pi \|_f$ for some appropriately defined function $f$.

Because the chain is irreducible and positive recurrent, $\tau_0 < \infty$ $\PP_x$-\as\, for $x \in \Nset$. By construction, for all
$x =1,2,\dots$,
$\tau_{x-1} \leq \tau_0$, $\PP_x$-\as, which implies that $\PE_x[\tau_0]= \PE_x[\tau_{x-1}] + \PE_{x-1}[\tau_0]$ and,
for any $s \in \Cset$ such that $|s| \leq 1$, $\PE_x[s^{\tau_0}]= \PE_x[s^{\tau_{x-1}}]\PE_{x-1}[s^{\tau_0}]$, where $\tau_{x-1}$ is the
first return time of the state $x-1$. For all $x=1,2,\dots$, we have $\PP_x\{ \tau_{x-1} \in \cdot \}= \PP_1\{ \tau_0 \in \cdot\}$
which shows that $\PE_x[\tau_0]= x \PE_1[\tau_0]$ and $\PE_x[s^{\tau_0}] = e^x(s)$, where $e(s) \eqdef \PE_1[s^{\tau_0}]$.
This relation implies
\begin{equation*}
e(s) = s a_0 + \sum_{y=1}^\infty a_y e^y(s) 
= s \int_0^\infty \rme^{\lambda (e(s)-1)t} d B(t) \eqsp.
\end{equation*}
By differentiating the previous relation \wrt\ $s$ and taking the limit as $s \to 1$, the previous relation implies that:
$\PE_1[\tau_0]= (1-\rho)^{-1}$.
Since $\{0,1\}$ is an atom, we may use Theorem \ref{theo:mainresult:stochasticallymonotone} with $C= \{ 0,1\}$, $r \equiv 1$ and $\driftv_0 \equiv 1$.
In this case
\[
U_0(x) = V_0(x)= 1+\PE_x[\sigma_C]= 1+\PE_{x-1}[\tau_0] \1_{ \{ x \geq 2 \}} =1+(1-\rho)^{-1}(x-1) \1_{ \{ x \geq 2 \}} \eqsp.
\]
Theorem \ref{theo:mainresult} shows that, for any $(x,x') \in \Nset \times \Nset$ and any functions $\alpha$ and $\beta$ satisfying \eqref{eq:additivity},
\[
\alpha(n) \| P^n(x,\cdot) - P^n(x',\cdot) \|_{\beta} \leq 1+(1-\rho)^{-1}(x \vee x' -1) \1_{ \{ x \vee x' \geq 2 \}} \eqsp.
\]
Convergence bounds $\alpha(n) \| P^n(x,\cdot)-\pi \|_\beta$ can be obtained by integrating the previous relation in $x'$
\wrt\ the stationary distribution $\pi$ (which can be computed using the Pollaczek-Khinchine formula).

\begin{figure}
\centering
  \includegraphics[width=9cm]{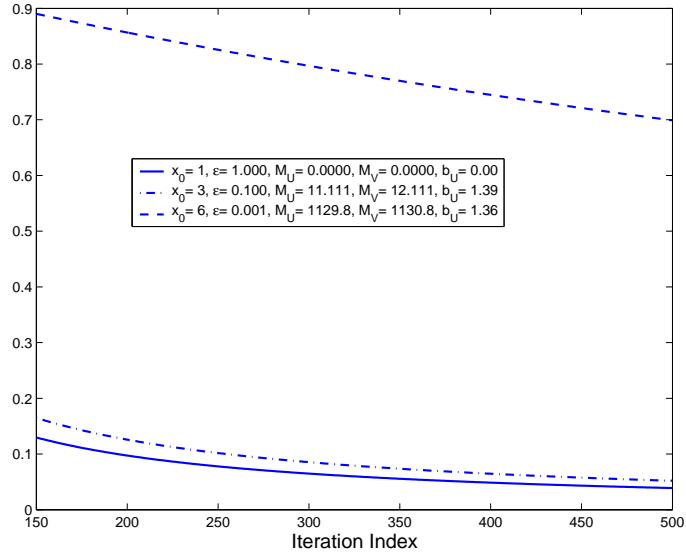}\\
  \caption{convergence bound for the total variation distance in the light-traffic case: $\rho=0.5$, $\alpha=2.5$  }
  \label{fig:total-variation-MG1-0.5}
\end{figure}

\begin{figure}
  \centering
  \includegraphics[width=9cm]{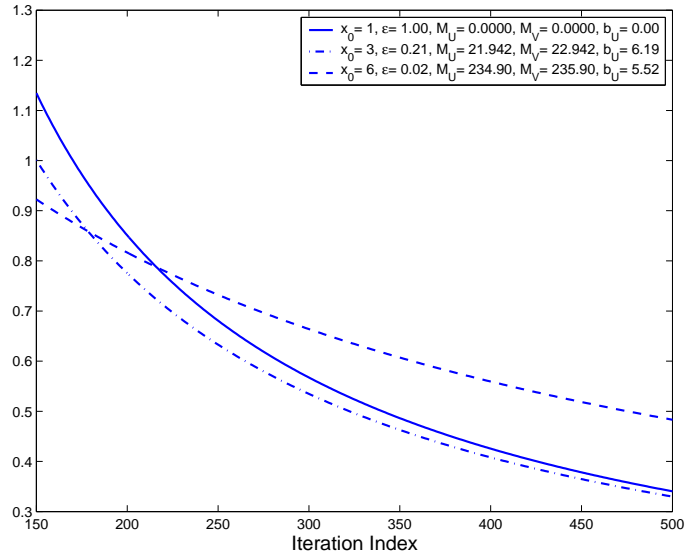}\\
  \caption{convergence bound for the total variation distance in the heavy traffic case: $\rho= 0.9$, $\alpha=2.5$}
  \label{fig:total-variation-MG1-0.9}
\end{figure}

It is possible to choose the set $C$ in a different way, leading to different bounds.  One may set for example $C= \{0, \dots, x_0\}$, for
some $x_0 \geq 2$.  For simplicity,
assume that the sequence $\{a_j \}_{j \geq 0}$ is non-decreasing. In this case, for all $x \in C$ and $y \in \Nset$,
$P(x,y) = a_{y-x+1} \1_{\{y \geq x-1\}} \geq a_y \1_{\{ y \geq x_0-1 \}}$ and the set $C$ satisfies \ref{assumption:smallset} with $\epsilon \eqdef \sum_{y=x_0-1}^\infty a_y$ and
$\nu(y)= \epsilon^{-1} a_y \1_{\{ y \geq x_0-1\}}$. Taking again $r (k) \equiv 1$ and $\driftv_0(x) \equiv 1$, we have
\begin{multline*}
U_0(x) = V_0(x) = 1 + \PE_x[\tau_C] \1_{C^c}(x) = 1 + \PE_x[\tau_{x_0}] \1_{C^c}(x) \\
= 1 + \PE_{x-x_0}[\tau_0] \1_{C^c}(x) = 1 +  (1-\rho)^{-1} (x-x_0) \1_{C^c}(x) \eqsp.
\end{multline*}
To apply the results of Theorems \ref{theo:mainresult:stochasticallymonotone}, we finally compute a  bound for
$b_{U_0}= \sup_C Q U_0 = (1-\epsilon)^{-1} [\sup_C P U_0 - \epsilon \nu(U_0)]$, which can be obtained by combining a bound
for $\sup_C P U_0$ and the expression of $\nu(U_0)$. An expression $\nu(U_0)$ is computed by a direct application of the definitions.
The bound for $\sup_C P U_0$ is obtained by noting that, for all $y > x_0$ and $x \in C$, $P(x,y) \leq P(x_0,y)= a_{y-x_0+1}$, which implies
\begin{multline*}
P U_0(x) = \PE_x[\tau_C]= 1 + \PE_x \left[ \PE_{X_1}[\tau_C] \1_{\{\tau_C > 1\}}\right]  = 1 + \PE_x \left[ \PE_{X_1}[\tau_{x_0}]  \1_{\{X_1 \not \in C \}} \right] \\
= 1 + (1-\rho)^{-1} \sum_{y=x_0+1}^\infty (y-x_0) P(x,y) \leq 1 + (1-\rho)^{-1} \sum_{y=x_0+1}^\infty (y-x_0) a_{y-x_0+1} \eqsp.
\end{multline*}
We provide some numerical illustrations of the bounds described above. We use the distribution of service time suggested by
in \cite{roughan:veitch:rumsewicz:1998} given by
\begin{equation}
\label{eq:servicetime}
b(x)=
\begin{cases}
\alpha B^{-1} \rme^{-\frac{\alpha}{B}x} & x \leq B \\
\alpha B^\alpha \rme^{-\alpha} x^{-\alpha+1} & x > B
\end{cases}
\end{equation}
where $B$ marks where the tail begins. The mean of the  service distribution is $m_1= B \left\{1 + \rme^{-\alpha}/(\alpha-1) \right \}/\alpha$
and its Laplace transform, $G(s) = \int_0^\infty \rme^{-st} d B(t)$, $s \in \Cset$, $\mathrm{Re}(s) \geq 0$, is given by
\[
G(s)= \alpha \frac{1-\rme^{-(sB+\alpha)}}{sB+\alpha} + \alpha B^\alpha \re^{-\alpha} s^\alpha \Gamma(-\alpha,sB) \eqsp,
\]
where $\Gamma(x,z)$ is the incomplete $\Gamma$ function. The probability generating function $P_\pi(z)$
of the stationary distribution is given by the Pollaczek-Khinchine formula
\[
P(z) = \frac{(1-\rho)(z-1) G(\lambda(1-z))}{z-G(\lambda(1-z))} \eqsp.
\]
In figures \ref{fig:total-variation-MG1-0.5} and \ref{fig:total-variation-MG1-0.9},
we display the convergence bound $\|P^n(x,\cdot) - \pi \|_{\mathrm{TV}}$ as a function of the iteration index $n$,
for $x=10$, $\alpha= 2.5$, different choices of the small set upper limit $x_0=1,3,6$, and two different values of the traffic $\rho= 0.5$ (light traffic)
and $\rho=0.9$ (heavy traffic). Perhaps surprisingly, the bound computed using the atom $C= \{0,1\}$ is not better uniformly in the iteration index
$n$. There is a trade off between the number of visits to the small set where coupling might and the probability that coupling is successful.
In the heavy traffic case ($\rho= 0.9$), the queue is not  very often empty, so the atom is not frequently visited, explaining why deriving the
convergence bound from a larger coupling set improves the bound (this effect is even more noticeable for a critically loaded system).
\begin{center}
Insert figures \ref{fig:total-variation-MG1-0.5} and \ref{fig:total-variation-MG1-0.9} approximately here
\end{center}

\subsection{The Independence Sampler}
\label{subsec:independencesampler}
This second example is borrowed from \cite{jarner:roberts:2001}. It is
an example of a Markov chain which is stochastically
monotone w.r.t a non-standard ordering of the state and does
not have an atom at the bottom of the state-space.

The purpose of the Metropolis-Hastings Independence Sampler is to
sample from a probability density $\pi$ (\wrt\ some $\sigma$-finite measure
$\mu$ on $\Xset$), which is known only up to a scale factor.  At each iteration, a move is proposed according to a
distribution with density $q$ with respect to $\mu$.  The move
is accepted with probability $a(x,y) \eqdef \frac{q(x)}{\pi(x)} \frac{\pi(y)}{q(y)} \; \wedge 1$.
The transition kernel of the algorithm is thus given by
\[
P(x,A) = \int_A a(x,y) q(y) \, \mu(dy) + \1_{A}(x) \int_\Xset
\Big(1 -a(x,y) \Big) q(y) \, \mu(dy), \quad x \in \Xset, A \in
\Xsigma.
\]
It is well known that the independence sampler is stochastically monotone with respect to the ordering:
$ x' \preceq x \Leftrightarrow \frac{q(x)}{\pi(x)} \leq \frac{q(x')}{\pi(x')}$.
Without
loss of generality, it is assumed that $\pi(x) > 0$ for all $x \in
\Xset$ and  that $q > 0$ $\pi$-\as .  For all $\eta>0$, define the set
\begin{equation}
\label{eq:definitionCepsilon}
C_\eta \eqdef \left\{x \in \Xset: \frac{q(x)}{\pi(x)}\geq\eta \right\} \eqsp.
\end{equation}
For any $\eta > 0$, we assume that $0 < \pi(C_\eta) < 1$ and we denote by
$\nu_\eta(\cdot)$ the probability measure $\nu_\eta(\cdot)= \pi(\cdot \cap C_\eta)/ \pi(C_\eta)$.
For any $x\in C_\eta$,
\begin{multline*}
  P(x,A) \geq \int_A \left( \frac{q(x)}{\pi(x)} \wedge
    \frac{q(y)}{\pi(y)} \right) \pi(y) \mu(dy) \\
  \geq \int_{A \cap C_\eta} \left( \frac{q(x)}{\pi(x)} \wedge
    \frac{q(y)}{\pi(y)} \right) \pi(y) \mu(dy) \geq \eta \pi(A
  \cap C_\eta) = \eta \pi(C_\eta) \nu_\eta(A).
\end{multline*}
showing that the set $C_\eta$ satisfies \ref{assumption:smallset} with $\nu = \nu_\eta$ and $\epsilon= \eta \pi(C_\eta)$.
\begin{prop}  \label{prop:withouttears}
Assume that there exists a  decreasing differentiable function $K: (0,\infty) \to (1,\infty)$,
whose inverse is denoted by $K^{-1}$, satisfying
\begin{enumerate}
\item \label{item:condIS1} the function $\phi(v)= v K^{-1}(v)$ is differentiable, increasing
  and concave on $[1,\infty)$, $\lim_{v \to \infty} \phi(v) = \infty$, and $\lim_{v \to \infty} \phi'(v) = 0$.
\item \label{item:condIS2} $\int_0^{+\infty}  u K(u) d \psi(u) < \infty$, where for $\eta>0$, $\psi(\eta) \eqdef 1-\pi(C_\eta)$.
\end{enumerate}
Then, for any $\eta^\star $ satisfying
\[
 \left\{ 1-\psi(\eta^\star) \right\} \phi(1) > \int_0^\infty (u \wedge \eta^\star) K(u) d \psi(u)
\]
assumption \ref{assumption:driftfunction-stochasticallymonotone} is satisfied with $\driftfunction_0= K \circ (q/\pi)$, $C=C_{\eta^\star}$ and
\[
\phi_0(v) = \{ 1-\psi(\eta^\star) \} \phi(v) - \int_0^\infty(u \wedge \eta^\star) K(u) d \psi(u) \eqsp.
\]
In addition,
\[
\sup_{x \in C_{\eta^\star}} P \driftfunction_0 \leq \int_0^{+\infty}  u K(u) d \psi(u)  + K(\eta^\star) \eqsp.
\]
\end{prop}
To illustrate our results, we evaluate the convergence bounds in the case where the target density  $\pi$ is the uniform distribution on $[0,1]$
and the proposal density is  $q(x)=(r+1)x^r \1_{[0,1]}(x)$.
Proposition \ref{prop:withouttears} provides a mean to derive a drift condition of the form $PW_0\leq W_0-\phi \circ W_0$ outside some small set $C$
for functions $\phi \in \mathcal{C}$ of the form $\phi(v)=c v^{1-1/\alpha}+d$ for any $\alpha \in  [1,1+1/r)$.
In this case, the function $\psi$ is given by $\psi(\eta)=(\eta/(r+1))^{1/r}$, for $\eta \in [0,r+1]$ and $\psi(\eta)= 1$ otherwise.
We set, for $u \in [0,r+1]$, $K(u)=(u/(r+1))^{-\alpha}$. The integral
$\int u K(u) d\psi(u)=\frac{(r+1)^{-\alpha}}{r(-\alpha+1/r+1)}$ is finite provided that $\alpha< 1+1/r$. The function
$\phi(u)=uK^{-1}(u)=u^{1-1/\alpha}(r+1)$ belongs to $\mathcal{C}$ provided that $\alpha>1$.

\begin{figure}
  \includegraphics[width=11cm]{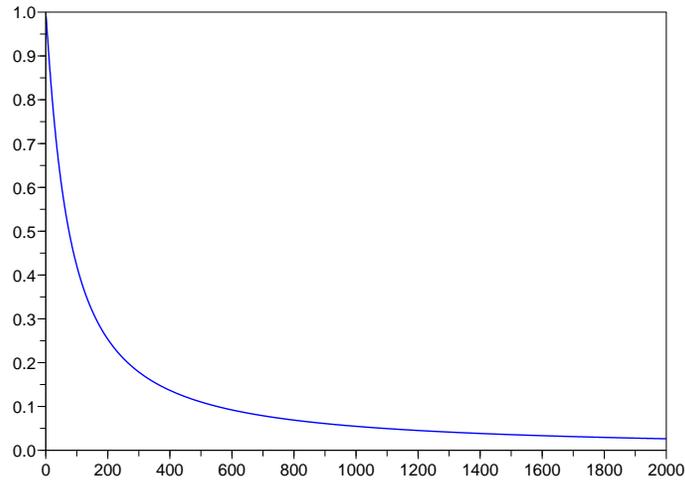}\\
  \caption{convergence bound for the total variation distance for the independence sampler with $q(x)=3x^2$}
  \label{fig:total-variation-IS-2}
\end{figure}

\begin{figure}
  \includegraphics[width=11cm]{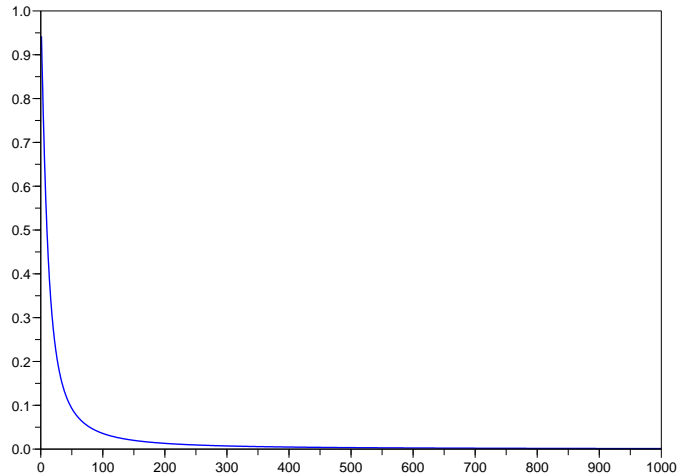}\\
  \caption{convergence bound for the total variation distance when $q(x)=1.5\sqrt{x}$}
  \label{fig:total-variation-IS-1/2}
\end{figure}

Using these results, it is now straightforward to evaluate the
constants in Theorem \ref{theo:mainresult}; this can be employed to
calculate a bound on exactly how many iterations are necessary to
get within a prespecified total variation distance of the target
distribution. In figures \ref{fig:total-variation-IS-2}  and
\ref{fig:total-variation-IS-1/2}, we have displayed the total
variation bounds to convergence for the instrumental densities $q(x) = 3x^2$ ($r=2$) and
$q(x)= (3/2) \sqrt{x}$. We have taken $\alpha=1.1$ and $\eta^\star= 0.25$ for $r=2$ and $\alpha=1.5$ and $\eta^\star= 0.5$ for $r=1/2$.
 When ($r=2$, $\alpha=1.1)$ the convergence to stationarity is
quite slow, which is not surprising since the instrumental density
does not match well the target density at zero $x=0$: according to our
computable bounds, $500$ iterations are required to get the total variation to
the stationary distribution below $0.1$. When $r=1/2$, the degeneracy of the
instrumental density at zero is milder and the convergence rate is significantly faster.
Less than $50$ iterations are required to reach the same bound.

\section{Proof of Theorem \ref{theo:mainresult}}
\label{sec:proofmainresult}
The proof is based on the pathwise coupling construction.
For $(x,x') \in \Xset \times \Xset$, and $A \in \Xsigma \otimes \Xsigma$, define $\bar{P}$ the coupling kernel as follows
\begin{align*}
& \bar{P}\left(x,x',0;A \times \{ 0 \} \right) = \left(1 - \epsilon \1_{C \times C}(x,x') \right) \Pdouble(x,x',A) \\
& \bar{P}\left(x,x',0;A \times \{ 1 \} \right) = \epsilon \1_{C \times C}(x,x') \nu(A \cap \{ (x,x') \in \Xset \times \Xset, x=x'\}) \\
& \bar{P}\left(x,x',1;A \times \{ 0 \} \right) = 0 \\
& \bar{P}\left(x,x',1;A \times \{ 1 \} \right) = \int P(x,dy)  \1_A(y,y) \eqsp.
\end{align*}
For any probability measure $(x,x') \in \Xset \times \Xset$,
denote $\bar{\PP}_{x,x'}$ and $\bar\PE_{x,x'}$ the probability measure and the expectation on associated to the Markov chain
$\{(X_n,X'_n,d_n)\}_{n \geq 0}$ with transition kernel $\bar{P}$ starting from $(X_0,X'_0,0) = (x,x',0)$.
In words, the coupling construction proceeds as follows. If $d_n=0$ and $(X_n,X'_n) \not \in C \times C$, we draw $(X_{n+1},X'_{n+1})$
according to $\Pdouble(x,x',\cdot)$ and set $d_{n+1}=0$. If $d_n=0$
and $(X_n,X'_n) \in C \times C$, we draw a coin with probability of heads $\epsilon$. If the coin comes up head, then we draw $X_{n+1}$
from $\nu$ and set $X'_{n+1}=X_{n+1}$ and $d_{n+1}=1$ (the coupling is said to be successful); if the coin comes up tails, then we draw $(X_{n+1},X'_{n+1})$
from $\Pdouble(X_n,X'_n,\cdot)$ and we set $d_{n+1}=1$. Finally, if $d_n=1$, we draw $X_{n+1}$ from $P(X_n,\cdot)$ and set $X_{n+1}=X'_{n+1}$.

By construction, for any $n$, $(x,x') \in \Xset \times \Xset$ and  $(A,A') \in \Xsigma \times \Xsigma$,
\begin{multline*}
  \bar\PP_{x,x',0}(Z_n \in A \times \Xset \times \{ 0,1 \}) = \bar\PP_{x,x',0}( X_n \in A) = P^n(x,A) \quad \text{and} \\
  \bar\PP_{x,x',0}(Z_n \in \Xset \times A' \times \{0,1\})= \bar\PP_{x,x',0}(X'_n \in A') = P^n(x',A')\eqsp.
\end{multline*}
By \cite[Lemma 1]{douc:moulines:rosenthal:2004}, we may relate the expectations of functionals
under the two probability measures $\bar\PP_{x,x',0}$ and $\PPdouble_{x,x'}$, where $\PPdouble_{x,x'}$ is defined in \eqref{eq:defPdouble}:
for any non-negative adapted process $(\chi_k)_{k\geq 0}$ and $(x,x') \in \Xset \times \Xset$,
\begin{gather} \label{eq:stop}
  \bar \PE_{x,x',0} [ \chi_{n} \1_{ \{T > n \} } ] = \PEdouble_{x,x'} \left[\chi_{n} \, (1-\epsilon)^{N_{n-1}} \right] \eqsp,
\end{gather}
where $N_n$ is the number of visit to the set ${C \times C}$ before time $n$,
\begin{equation}
\label{eq:definitionNk}
N_n = \sum_{j=0}^\infty \1_{\{ \sigma_j \leq n\}} = \sum_{i=0}^n \1_{C \times C}(X_i,X'_i) \eqsp.
\end{equation}
Let $f  : \Xset \to [0,\infty)$  and let $g : \Xset \to \Rset$ be any Borel function such that
$\sup_{x \in \Xset} |g(x)| / f(x) < \infty$. The classical coupling inequality (see e.g.\
\cite[Chapter 2, section 3]{thorisson:2000})
implies that
\begin{multline*}
\left | P^n(x,g) - P^n(x',g)  \right |
= \left | \bar \PE_{x,x',0} \left[ g(X_n) - g(X'_n) \right]  \right|
\\ \leq \sup_{x \in \Xset} |g(x)| / f(x) \; \bar \PE_{x,x',0} \left[ (f(X_n) + f(X'_n)) \1 \{ d_n = 0 \} \right] \eqsp,
\end{multline*}
and \eqref{eq:stop} shows the following key coupling inequality:
\begin{equation}
\label{eq:keycouplingiinequality}%
\left \| P^n(x,\cdot) - P^n(x',\cdot) \right\|_f \leq
 \; \PEdouble_{x,x'} \left \{ (f(X_n) + f(X'_n)) (1-\epsilon)^{N_{n-1}} \right \} \eqsp.
\end{equation}

%
%
%
Because by definition $\alpha(u) \beta(v) \leq \rho u+ (1-\rho) v$ for all $(u,v) \in \Rset^+ \times \Rset^+$ and
any non negative function $f$ satisfying
$f(x) + f(x') \leq \beta \circ \driftV(x,x')$
for all $(x,x') \in \Xset \times \Xset$, the coupling inequality \eqref{eq:keycouplingiinequality} shows that
\begin{align*}
  &\alpha \circ \left\{ R(n) + M_{\driftU} \right\}  \| P^n(x,\cdot) - P^n(x',\cdot) \|_{f} \\
  & \quad \leq \alpha \circ \left\{ R(n) + M_{\driftU} \right\} \PEdouble_{x,x'} [ \{f(X_n) + f(X'_n)\} (1-\epsilon)^{N_{n-1}}] \\
  & \quad \leq \rho \left\{ R(n) + M_{\driftU} \right\}\; \PEdouble_{x,x'} [    (1-\epsilon)^{N_{n-1}}] + (1-\rho) \; \PEdouble_{x,x'} [ \driftV(X_n,X'_n) (1-\epsilon)^{N_{n-1}}] \eqsp.
\end{align*}
Set for any $n \geq 0$, $\driftU_n(x,x')= \PEdouble_{x,x'} \left[ \sum_{k=0}^{\sigma_{C \times C}} r(n+k) \right]$.
It is well known that $\{ \driftU_n \}_{n \geq 0}$ satisfies the sequence of drift equations
\begin{equation}
\label{eq:driftequationVn}
\Pdouble \driftU_{n+1} \leq \driftU_n - r(n) + b_{\driftU} r(n) \1_{C \times C} \eqsp,
\end{equation}
Similarly, $\Pdouble \driftV  \leq \driftV  - \driftv + b_{\driftV} \1_{C \times C}$.
Define  for $n\geq 0$,
\begin{align*}
W^{(0)}_n   &\eqdef \driftU_n(X_n,X'_n) + \sum_{k=0}^{n-1} r(k) + M_{\driftU} \eqsp,\\
W^{(1)}_n &\eqdef \driftV(X_n,X'_n) + \sum_{k=0}^{n-1}  \driftv(X_k,X'_k) + M_{\driftV} \eqsp.
\end{align*}
with the convention $\sum_u^v = 0$ when $u > v$.

Since by construction, for any $n \geq 1$, $W^{(0)}_n \geq R(n)$ and $W^{(1)}_n \geq  \driftV(X_n,X'_n)$, the previous inequality implies,
\begin{multline*}
  \alpha \circ R(n) \| P^n(x,\cdot) - P^n(x',\cdot) \|_{f} \\
  \leq \PE_{x,x'}^\star [W^{(0)}_n (1-\epsilon)^{N_{n-1}}] + \PE_{x,x'}^\star [W^{(1)}_n (1-\epsilon)^{N_{n-1}}] \eqsp.
\end{multline*}
We now have to compute bounds for $\PE_{x,x'}^\star [W^{(i)}_n (1-\epsilon)^{N_{n-1}}]$, $i=0,1$.  Define
\begin{equation}
T^{(0)}_n \eqdef  \prod_{i=0}^{n-1} \frac{ W^{(0)}_{i} + b_{\driftU} r(i) \1_{ {C \times C}}( X_{i},X'_i)}{ W^{(0)}_{i} }  \quad \text{and} \quad
T^{(1)}_n \eqdef  \prod_{i=0}^{n-1} \frac{ W^{(1)}_{i} + b_{\driftV}  \1_{ {C \times C}}( X_{i},X'_i)}{ W^{(1)}_{i} } \eqsp.
\label{eq:defsn}
\end{equation}
If $\epsilon=1$, $(1-\epsilon)^{N_{n-1}} = \1_{ \{ \sigma_0 \geq n \}}$, where $\sigma_0 = \inf\{n \geq 0 \mid (X_n,X'_n) \in {C \times C}\}$
is the first hitting time of the set ${C \times C}$: $T^{(i)}_n \1_{\{\sigma_0 \geq n\}} =  \1_{\{ \sigma_0 \geq n \}} \leq 1$.
Consider now the case $\epsilon<1$. By construction,  for $N_{n-1}= 0$, $T_n^{(i)}=1$ and for $N_{n-1} > 0$,
\begin{gather} \label{eq:expressiontn}
T^{(0)}_n =  \prod_{i=0}^{N_{n-1}-1} \frac{ W^{(0)}_{{\sigma}_i} + b_{\driftU} r(\sigma_i)} {W^{(0)}_{{\sigma}_i}}
\quad \text{and} \quad T^{(1)}_n =  \prod_{i=0}^{N_{n-1}-1} \frac{ W^{(1)}_{{\sigma}_i} + b_{\driftV}} {W^{(1)}_{{\sigma}_i}}
\end{gather}
where  $\sigma_i$ are the successive hitting time of the set ${C \times C}$ recursively defined by $\sigma_{j+1} = \inf\{n > \sigma_j \mid (X_n,X'_n) \in  {C \times C}\}$.
Because $W_n^{(0)} \geq R(n+1) + M_\driftU$, and $1 + b_\driftU r(n) / \{ R(n+1) + M_\driftU\} \leq 1/(1-\epsilon)$, for $N_{n-1} > 0 $, we have
\begin{equation}
\label{eq:boundB0}
  T^{(0)}_n (1-\epsilon)^{N_{n-1}} \leq \prod_{i=0}^{N_{n-1}-1} \left(\left\{ 1 + \frac{ b_{\driftU} \, r(\sigma_i) }{R(\sigma_i+1)+M_{\driftU}}
  \right\} (1-\epsilon)\right)
  \leq 1 \eqsp.
\end{equation}
Similarly, because  $W_n^{(1)} \geq M_\driftV$ and $1 + b_\driftV/ M_\driftV \leq 1/(1-\epsilon)$, we have $T^{(1)}_n (1-\epsilon)^{N_{n-1}} \leq 1$.
These two relations imply, for $i=0,1$,
\begin{align*}
&{\PE}^\star_{x,x'} \left[ W^{(0)}_n (1-\epsilon)^{N_{n-1}}  \right]  \leq {\PE}^\star_{x,x'} \left[ W^{(0)}_n \{T^{(0)}_n \}^{-1}  \right] \eqsp, \\
&{\PE}^\star_{x,x'} \left[ W^{(1)}_n (1-\epsilon)^{N_{n-1}}  \right]  \leq  {\PE}^\star_{x,x'} \left[ W^{(1)}_n \{T^{(1)}_n \}^{-1}  \right] \eqsp.
\end{align*}
It remains now to compute a bound for ${\PE}^\star_{x,x'} \left[ W^{(i)}_n \{T^{(i)}_n\}^{-1}  \right]$. By construction, we have for $n \geq 1$,
\begin{multline}
\label{eq:basic}
\PEdouble_{x,x'} \left [ W^{(0)}_n \{T_{n}^{(0)}\}^{-1} \mid \mathcal{F}_{n-1} \right ] \\
= \PEdouble_{x,x'} \left [ W^{(0)}_n \mid \mathcal{F}_{n-1} \right ] \frac{W^{(0)}_{n-1}}{W^{(0)}_{n-1} + b_{\driftU} r(n-1)
\1_{C \times C}(X_{n-1},X'_{n-1})} \{T_{n-1}^{(0)}\}^{-1},
\end{multline}
where $\mcf_n = \sigma \left\{  (X_0,X'_0),\dots, (X_n,X'_n) \right\}$.
Now, \eqref{eq:driftequationVn}  yield:
\begin{gather} \label{eq:driftw}
  \PEdouble_{x,x'} \left [ W^{(0)}_n \mid \mcf_{n-1} \right ] \leq W^{(0)}_{n-1} + b_\driftU r(n-1) \1_{C \times C}(X_{n-1},X'_{n-1}) \eqsp.
\end{gather}
Combining \eqref{eq:basic} and \eqref{eq:driftw} shows that $\left \{W^{(0)}_n \{T_n^{(0)}\}^{-1} \right\}_{n \geq 0}$ is a $\mathcal F$-supermartingale.
Thus,
\begin{align*}
  {\PE}^\star_{x,x'} \left[ W^{(0)}_n (1-\epsilon)^{N_{n-1}} \right] \leq {\PE}^\star_{x,x'} \left[ W^{(0)}_n \{T_n^{(0)}\}^{-1} \right]
  \leq \PE_{x,x'}^\star[W^{(0)}_0] =  \driftU_0(x,x') + M_\driftU \eqsp.
\end{align*}
Similarly, ${\PE}^\star_{x,x'} \left[ W^{(1)}_n (1-\epsilon)^{N_{n-1}} \right] \leq  \driftV(x,x') + M_\driftV$, which concludes the proof of Theorem
\ref{theo:mainresult}.

\section{Proof of Proposition \ref{prop:mainresultwithdrift}, Theorem \ref{theo:singledriftcondition}}
\label{sec:proofmainresults-auxiliary}
\begin{proof}[Proof of Proposition \ref{prop:mainresultwithdrift}]
By applying the comparison Theorem \cite[]{meyn:tweedie:1993} and \cite[Proposition
  2.2]{douc:fort:moulines:soulier:2004}, we obtain the following inequalities. Then, for all
  $(x,x') \in \Xset \times \Xset$,
\begin{align}
& \PEdouble_{x,x'} \left [ \sum_{k=0}^{\tau_{C \times C}-1} \rphiunnorm(k) \right ] \leq \driftfunction(x,x') - 1 + b \frac{\rphiunnorm(1)}{\rphiunnorm(0)} \1_{C \times C}(x,x') \eqsp, \label{eq:dfms2} \\
\label{eq:dfms1} & \PEdouble_{x,x'} \left [ \sum_{k=0}^{\tau_{C \times C}-1} \phi \circ \driftfunction(X_k,X'_k) \right ] \leq \driftfunction(x,x') + b \1_{C \times C}(x,x') \eqsp.
\end{align}
The sequence $\{ \rphiunnorm(k) \}_{k \geq 0}$ is log-concave. Therefore, for any $k\geq 0$, $\rphiunnorm(k+1) / \rphiunnorm(k) \leq \rphiunnorm(1) / \rphiunnorm(0)$.
Then, applying \eqref{eq:dfms2}, we obtain:
\begin{multline*}
  \PEdouble_{x,x'} \left[ \sum_{k=0}^{\sigma_{C \times C}} \rphiunnorm(k) \right]  = \rphiunnorm(0) + \PEdouble_{x,x'} \left[ \sum_{k=1}^{\tau_{C \times C}} \rphiunnorm(k) \right]\1_{(C \times C)^c}(x,x') \\ \leq \rphiunnorm(0) + \frac{\rphiunnorm(1)}{\rphiunnorm(0)}
  \, \PEdouble_{x,x'} \left[\sum_{k=1}^{\tau_{C \times C}} \rphiunnorm(k-1) \right ] \1_{(C \times C)^c}(x,x') \eqsp,
\end{multline*}
showing \eqref{eq:boundV0}. Similarly,
\begin{multline*}
  \PEdouble_{x,x'} \left[ \sum_{k=0}^{\sigma_{C \times C}} \phi \circ \driftfunction(X_k,X'_k) \right]   = \phi \circ \driftfunction(x,x') \1_{{C \times C}}(x,x') \\
  + \PEdouble_{x,x'} \left [\sum_{k=0}^{\tau_{C \times C}-1} \phi \circ \driftfunction(X_k,X'_k) \right] \1_{(C \times C)^c}(x,x')  +
  \PEdouble_{x,x'}[\phi\circ \driftfunction(X_\tau,X'_\tau)] \1_{(C \times C)^c}(x,x')  
\end{multline*}
showing \eqref{eq:boundV1}. 
\end{proof}

\begin{proof}[Proof of Theorem \ref{theo:singledriftcondition}]
Since $d_0 = \inf_{x \not \in C} W_0(x)$, if $(x,x') \not \in C\times C$, $\driftfunction(x,x') \geq d_0$ and $\1_C(x) + \1_C(x') \leq 1$ since
either $x \not \in C$, $x' \not \in C$ (or both). The definition of the
kernel $\Pdouble$ therefore implies
\begin{align*}
\Pdouble \driftfunction(x,x')
&\leq \driftfunction_0(x) + \driftfunction_0(x') - 1 - \phi_0 \circ \driftfunction_0(x') - \phi_0 \circ \driftfunction_0(x') + b_0 \left\{ \1_C(x) + \1_C(x') \right\} \\
&\leq \driftfunction(x,x') - \phi_0 \circ \driftfunction (x,x') + b_0  \eqsp,
\end{align*}
where we have used the inequality: for any $u \geq 1$ and $v \geq 1$, $\phi_0(u+v-1) - \phi_0(u) \leq \phi_0(v) - \phi_0(1)$. For $(x,x') \not \in C$,
$ b_0 \leq (1-\lambda) \phi_0(d) \leq (1-\lambda) \phi_0 \circ \driftfunction_0(x,x')$  and the previous inequality implies
$\Pdouble \driftfunction(x,x') \leq \driftfunction(x,x') - \phi \circ \driftfunction(x,x')$.
\end{proof}

\appendix
\section{Proof of Proposition \ref{prop:withouttears}}
Let $W$ be any measurable non negative function on $\Xset$.  Then, for
$\eta>0$ and $x \not \in C_\eta$,
\begin{align*}
  PW(x)-W(x) & =  \int_\Xset a(x,y) \{W(y)-W(x)\} q(y) \mu(dy) \\
  & \leq \int_\Xset \left( \eta \wedge
    \frac{q(y)}{\pi(y)}\right)W(y)\pi(y)\mu(dy) - W(x) \int_\Xset
  a(x,y) q(y) \mu(dy).
\end{align*}
If $x\notin C_\eta$ and $y \in C_\eta$, then $y \preceq x$ and
$a(x,y) q(y) = (q(x)/\pi(x)) \, \pi(y)$.  Thus, we have:
\[
\int_\Xset a(x,y) q(y) \mu(dy) \geq \int_{C_\eta} a(x,y)
q(y) \mu(dy) = \frac{q(x)}{\pi(x)} \pi(C_\eta)= \frac{q(x)}{\pi(x)} (1- \psi(\eta)).
\]
Altogether, we obtain, for all $x\notin C_\eta$:
\begin{align} \label{eq:drift1}
  PW(x) - W(x) \leq  \int_\Xset \left( \eta \wedge
    \frac{q(y)}{\pi(y)}\right)W(y)\pi(y)\mu(dy) - \{1 - \psi(\eta) \} \frac{q(x)}{\pi(x)} W(x).
\end{align}
Applying the definition of $W_0$, we now have:
\begin{multline} \label{eq:pivo}
  \int_\Xset \left( \eta \wedge
    \frac{q(y)}{\pi(y)}\right)W_0(y)\pi(y)\mu(dy) \\= \int_{\Xset} \left(
    \eta \wedge \frac{q(y)}{\pi(y)} \right)
  K\left(\frac{q(y)}{\pi(y)}\right) \pi(y) \mu(dy) = \int_0^\infty
  (\eta \wedge u) K(u) d\psi(u) < \infty.
\end{multline}
By Lebesgue's bounded convergence theorem, $\lim_{\eta\to0} \int_0^\infty
(\eta \wedge u) K(u) d\psi(u) = 0$. Since moreover $\lim_{\eta\to0}
\psi(\eta) = 0$, hence, for $\eta$ small enough, $\{1 -
\psi(\eta)\} \phi(M) > \int_0^\infty (\eta \wedge u) K(u) d\psi(u)$, hence
$\eta^\star$ is well defined. Now, \eqref{eq:drift1} and \eqref{eq:pivo} yield, for all $x \not \in
C_{\eta^\star}$,
\begin{align*}
  PW_0(x) - W_0(x) & \leq \int_0^\infty ( \eta^\star \wedge u) K(u)
  d\psi(u) - ( 1 - \psi(\eta^\star)) W_0(x) K^{-1} \circ W_0(x) \\
& = - \phi_0(W_0(x)).
\end{align*}
For $x \in C_{\eta^\star}$, we have $W_0(x) \leq K(\eta^\star)$.
Finally, we have, for any $x \in C_{\eta^\star}$,
\begin{align*}
  P W_0(x)  & \leq \int_\Xset q(y) W_0(y) \mu(dy) + W_0(x) \\
  & = \int_{\Xset} \frac{q(y)}{\pi(y)} K \left( \frac{q(y)}{\pi(y)}
  \right) \, \pi(y) \mu(dy) \, + \, W_0(x) \leq \int_0^\infty u K(u) d
  \psi(u) + K(\eta^\star).
\end{align*}

\bibliographystyle{ims}

\begin{thebibliography}{27}
\expandafter\ifx\csname natexlab\endcsname\relax\def\natexlab#1{#1}\fi
\expandafter\ifx\csname url\endcsname\relax
  \def\url#1{\texttt{#1}}\fi
\expandafter\ifx\csname urlprefix\endcsname\relax\def\urlprefix{URL }\fi

\bibitem[{Baxendale(2005)}]{baxendale:2005}
\textsc{Baxendale, P.~H.} (2005).
\newblock Renewal theory and computable convergence rates for geometrically
  ergodic {M}arkov chains.
\newblock \textit{Ann. Appl. Probab.} \textbf{15} 700--738.

\bibitem[{Douc et~al.(2004{\natexlab{a}})Douc, Fort, Moulines and
  Soulier}]{douc:fort:moulines:soulier:2004}
\textsc{Douc, R.}, \textsc{Fort, G.}, \textsc{Moulines, E.} and
  \textsc{Soulier, P.} (2004{\natexlab{a}}).
\newblock Practical drift conditions for subgeometric rates of convergence.
\newblock \textit{Ann. Appl. Probab.} \textbf{14} 1353--1377.

\bibitem[{Douc et~al.(2004{\natexlab{b}})Douc, Moulines and
  Rosenthal}]{douc:moulines:rosenthal:2004}
\textsc{Douc, R.}, \textsc{Moulines, E.} and \textsc{Rosenthal, J.}
  (2004{\natexlab{b}}).
\newblock Quantitative bounds for geometric convergence rates of {M}arkov
  chains.
\newblock \textit{Annals of Applied Probability} \textbf{14} 1643--1665.

\bibitem[{Fort(2001)}]{fort:2001}
\textsc{Fort, G.} (2001).
\newblock \textit{Contrôle explicite d'ergodicité de chaînes de Markov:
  applications à l'analyse de convergence de l'algorithme Monte-Carlo {EM}}.
\newblock Ph.D. thesis, Université de Paris VI.

\bibitem[{Fort and Moulines(2000)}]{fort:moulines:2000}
\textsc{Fort, G.} and \textsc{Moulines, E.} (2000).
\newblock {$V$}-subgeometric ergodicity for a {H}astings-{M}etropolis
  algorithm.
\newblock \textit{Statist. Probab. Lett.} \textbf{49} 401--410.

\bibitem[{Fort and Moulines(2003{\natexlab{a}})}]{fort:moulines:2003}
\textsc{Fort, G.} and \textsc{Moulines, E.} (2003{\natexlab{a}}).
\newblock Convergence of the {M}onte {C}arlo expectation maximization for
  curved exponential families.
\newblock \textit{Ann. Statist.} \textbf{31} 1220--1259.

\bibitem[{Fort and Moulines(2003{\natexlab{b}})}]{fort:moulines:2003:SPA}
\textsc{Fort, G.} and \textsc{Moulines, E.} (2003{\natexlab{b}}).
\newblock Polynomial ergodicity of {M}arkov transition kernels,.
\newblock \textit{Stochastic Processes and Their Applications} \textbf{103}
  57--99.

\bibitem[{Gulinsky and Veretennikov(1993)}]{gulinsky:veretennikov:1993}
\textsc{Gulinsky, O.~V.} and \textsc{Veretennikov, A.~Y.} (1993).
\newblock \textit{Large deviations for discrete-time processes with averaging}.
\newblock VSP, Utrecht.

\bibitem[{Jarner and Roberts(2001)}]{jarner:roberts:2001}
\textsc{Jarner, S.} and \textsc{Roberts, G.~O.} (2001).
\newblock Polynomial convergence rates of {M}arkov chains.
\newblock \textit{Annals of Applied Probability} \textbf{12} 224--247.

\bibitem[{Klokov and Veretennikov(2004)}]{klokov:veretennikov:2004}
\textsc{Klokov, S.~A.} and \textsc{Veretennikov, A.~Y.} (2004).
\newblock Sub-exponential mixing rate for a class of {M}arkov chains.
\newblock \textit{Math. Commun.} \textbf{9} 9--26.

\bibitem[{Lindvall(1979)}]{lindvall:1979}
\textsc{Lindvall, T.} (1979).
\newblock On coupling of discrete renewal sequences.
\newblock \textit{Z. Wahrsch. Verw. Gebiete} \textbf{48} 57--70.

\bibitem[{Lindvall(1992)}]{lindvall:1992}
\textsc{Lindvall, T.} (1992).
\newblock \textit{Lectures on the Coupling Method}.
\newblock Wiley, New-York.

\bibitem[{Lund et~al.(1996)Lund, Meyn and Tweedie}]{lund:meyn:tweedie:1996}
\textsc{Lund, R.~B.}, \textsc{Meyn, S.~P.} and \textsc{Tweedie, R.} (1996).
\newblock Computable exponential convergence rates for stochastically ordered
  {M}arkov processes.
\newblock \textit{Annals of Applied Probability} \textbf{6} 218--237.

\bibitem[{Lund and Tweedie(1996)}]{lund:tweedie:1996}
\textsc{Lund, R.~B.} and \textsc{Tweedie, R.~L.} (1996).
\newblock Geometric convergence rates for stochastically ordered {M}arkov
  chains.
\newblock \textit{Mathematics of Operation Research} \textbf{21} 182--194.

\bibitem[{Meyn and Tweedie(1993)}]{meyn:tweedie:1993}
\textsc{Meyn, S.~P.} and \textsc{Tweedie, R.~L.} (1993).
\newblock \textit{{M}arkov Chains and Stochastic Stability}.
\newblock Springer, London.

\bibitem[{Meyn and Tweedie(1994)}]{meyn:tweedie:1994}
\textsc{Meyn, S.~P.} and \textsc{Tweedie, R.~L.} (1994).
\newblock Computable bounds for convergence rates of {M}arkov chains.
\newblock \textit{Annals of Applied Probability} \textbf{4} 981--1011.

\bibitem[{Nummelin and Tuominen(1983)}]{nummelin:tuominen:1983}
\textsc{Nummelin, E.} and \textsc{Tuominen, P.} (1983).
\newblock The rate of convergence in {O}rey's theorem for {H}arris recurrent
  {M}arkov chains with applications to renewal theory.
\newblock \textit{Stochastic Processes and Their Applications} \textbf{15}
  295--311.

\bibitem[{Roberts and Rosenthal(2004)}]{roberts:rosenthal:2004}
\textsc{Roberts, G.~O.} and \textsc{Rosenthal, J.~S.} (2004).
\newblock General state space {M}arkov chains and {MCMC} algorithms.
\newblock \textit{Probab. Surv.} \textbf{1} 20--71.

\bibitem[{Roberts and Tweedie(1999)}]{roberts:tweedie:1999}
\textsc{Roberts, G.~O.} and \textsc{Tweedie, R.~L.} (1999).
\newblock Bounds on regeneration times and convergence rates for {M}arkov
  chains.
\newblock \textit{Stochastic Processes and Their Applications} \textbf{80}
  211--229.

\bibitem[{Roberts and Tweedie(2000)}]{roberts:tweedie:2000}
\textsc{Roberts, G.~O.} and \textsc{Tweedie, R.~L.} (2000).
\newblock Rates of convergence of stochastically monotone and continuous time
  {M}arkov models.
\newblock \textit{Journal of Applied Probability} \textbf{37} 359--373.

\bibitem[{Rosenthal(1995)}]{rosenthal:1995}
\textsc{Rosenthal, J.~S.} (1995).
\newblock Minorization conditions and convergence rates for {M}arkov chain
  {M}onte {C}arlo.
\newblock \textit{J. Am. Statist. Assoc.} \textbf{90} 558--566.

\bibitem[{Roughan et~al.(1998)Roughan, Veitch and
  Rumsewicz}]{roughan:veitch:rumsewicz:1998}
\textsc{Roughan, M.}, \textsc{Veitch, D.} and \textsc{Rumsewicz, M.} (1998).
\newblock Computing queue-length distributions for power-law queues.
\newblock In \textit{Proceedings. IEEE INFOCOM '98}, vol.~1. IEEE.

\bibitem[{Scott and Tweedie(1996)}]{scott:tweedie:1996}
\textsc{Scott, D.~J.} and \textsc{Tweedie, R.~L.} (1996).
\newblock Explicit rates of convergence of stochastically ordered {M}arkov
  chains.
\newblock In \textit{Athens Conference on Applied Probability and Time Series:
  {A}pplied Probability in Honor of J. M. {G}ani}, vol. 114 of \textit{Lecture
  Notes in Statistics}. Springer.

\bibitem[{Thorisson(2000)}]{thorisson:2000}
\textsc{Thorisson, H.} (2000).
\newblock \textit{Coupling, Stationarity and Regeneration}.
\newblock Probability and its Applications, Springer-Verlag, New-York.

\bibitem[{Tuominen and Tweedie(1994)}]{tuominen:tweedie:1994}
\textsc{Tuominen, P.} and \textsc{Tweedie, R.} (1994).
\newblock Subgeometric rates of convergence of $f$-ergodic {M}arkov {C}hains.
\newblock \textit{Advances in Applied Probability} \textbf{26} 775--798.

\bibitem[{Veretennikov(1997)}]{veretennikov:1997}
\textsc{Veretennikov, A.} (1997).
\newblock On polynomial mixing bounds for stochastic differential equations.
\newblock \textit{Stochastic Process. Appl.} \textbf{70} 115--127.

\bibitem[{Veretennikov(1999)}]{veretennikov:1999}
\textsc{Veretennikov, A.} (1999).
\newblock On polynomial mixing and the rate of convergence for stochastic
  differential and difference equations.
\newblock \textit{Theory of probability and its applications}  361--374.

\end{thebibliography}

\newpage
\clearpage

\end{document}